%% file: paper_SD_QP.tex
\RequirePackage{fix-cm}

\documentclass[smallextended]{svjour3}       

\smartqed  

\usepackage{cite}
\usepackage{algpseudocode}
\usepackage[ruled]{algorithm}
 \usepackage{url}
\usepackage{booktabs}
\usepackage{graphicx}
\usepackage{adjustbox}
\usepackage{tabularx}
\usepackage{siunitx}
\usepackage{color}
\usepackage{latexsym}
\usepackage{amsfonts}
\usepackage{amsmath}
\usepackage{subfigure}
\usepackage{multirow}
\usepackage{mathtools}
\usepackage{enumerate}

\newcommand{\R}{\rm{I\!R}}

\newcommand{\eps}{\varepsilon}

\def\square{{\setbox0=\hbox{X}\hbox to \ht0{\vrule\hss\vbox to \ht0{
  \hrule width \ht0\vfil\hrule width \ht0}\vrule}}}

\makeatletter

\begin{document}

\title{A simplicial decomposition framework for large scale convex quadratic programming}

\titlerunning{A simplicial decomposition framework for convex quadratic programming}        

\author{Enrico Bettiol         \and
	Lucas L\'etocart         \and
	Francesco Rinaldi         \and
	        Emiliano Traversi 
}


\institute{E. Bettiol \at
              LIPN, CNRS, (UMR7030), Universit\'e Paris 13, Sorbonne Paris Cit\'e, \\
              99 av. J. B. Cl\'ement, 93430 Villetaneuse, France \\
              \email{bettiol@lipn.univ-paris13.fr}           
           \and
           L. L\'etocart \at
              LIPN, CNRS, (UMR7030), Universit\'e Paris 13, Sorbonne Paris Cit\'e, \\
              99 av. J. B. Cl\'ement, 93430 Villetaneuse, France \\
              \email{lucas.letocart@lipn.univ-paris13.fr}           
           \and
           F. Rinaldi \at
              Dipartimento di Matematica, Universit\`a di Padova \\
              Via Trieste, 63, 35121 Padua, Italy \\
              Tel.: +39-049-8271424\\
              \email{rinaldi@math.unipd.it}           
           \and
           E. Traversi \at
              LIPN, CNRS, (UMR7030), Universit\'e Paris 13, Sorbonne Paris Cit\'e, \\
              99 av. J. B. Cl\'ement, 93430 Villetaneuse, France \\
              \email{emiliano.traversi@lipn.univ-paris13.fr}
}

\date{Received: date / Accepted: date}

\maketitle

\begin{abstract}

In this paper, we analyze in depth a simplicial decomposition like algorithmic framework for large scale
convex quadratic programming. In particular, we first propose two tailored strategies
for handling the master problem. Then, we describe a few techniques for speeding up the solution of the pricing  
problem. We report extensive numerical experiments  on both real portfolio optimization and general quadratic programming problems
showing the efficiency and robustness of the method when compared to  \texttt{Cplex}.


\end{abstract}

\section{Introduction}\label{sec:introduction}

We consider the following problem 
\begin{equation}\label{prob}
\begin{array}{lll}
\min& f(x)~=~ x^{\top} Q x + c^{\top} x  \\
\mbox{ s.t. }&  a_i^Tx~=~b_i,&i\in \cal{E}\\
&  a_i^Tx~\ge~b_i,&i\in \cal{I}
\end{array}
\end{equation}
with $Q\in \R^{n \times n}$, $c\in \R^{n}$, $a_i\in \R^{n}$ and $b_i \in \R, \ i \in\cal{E}\cup\cal{I}$.

Moreover, we assume that the polyhedral set 
\begin{equation*}
 X=\{x\in {\R}^n:~a_i^Tx~=~b_i,~i\in {\cal{E}} \} \cap\{x\in {\R}^n:~~a_i^Tx~\ge~b_i,~i \in {\cal{I}}\} 
\end{equation*}
is non-empty and bounded and that the Hessian matrix $Q$ is positive semidefinite. 
Among all possible problems of type~\eqref{prob}, we are particularly interested in the ones with the following additional properties:
\begin{itemize}
 \item The number of constraints is  considerably  smaller than the number of variables in the problem, i.e. $|{\cal{E}}\cup{\cal{I}}|\ll n$;
 \item the Hessian matrix $Q$ is dense.
	\end{itemize}
A significant number of large-scale problems, arising in many different fields (e.g. Communications, Statistics, Economics and Machine Learning), present a structure similar to 
the one described above \cite{boyd2004convex}.

Solution methods for this class of problems can be mainly categorized into either interior point
methods  or  active  set  methods \cite{nocedal2006sequential}.
In interior point methods, a sequence of parametrized barrier functions is (approximately) minimized using Newton's method.  
The main computational burden is represented by the calculation of the Newton system solution (used to get the search direction). 
Even if those methods are relatively recent (they started becoming populare in the 1990s), a large number of papers and books exist
related to them (see, e.g. \cite{gondzio2012interior, wright2005interior, nesterov1994interior, wright1997primal, ye2011interior}).

In active set methods, at each iteration, a working set that estimates the set of active constraints at the solution
is iteratively updated. This gives a subset of constraints to watch while searching the solution (which obviously reduces the complexity 
of our search in the end). Those methods, which have been widely used since the 1970s, turn out to be effective when dealing with small- and medium-sized
problems. They usually guarantee efficient detection of unboundedness and infeasibility (other than returning an accurate estimate of the optimal active set).
An advantage of active set methods over interior points is that they are well-suited for warmstarts,
where a good estimate of the optimal active set or solution is used to initialize the algorithm.
This turns out to be extremely useful in applications where a sequence of QP problems is solved,
e.g., in a sequential quadratic programming method. A detailed overview of active set methods can be found in \cite{nocedal2006sequential}.

In this paper, we develop a \emph{simplicial decomposition type} approach (see e.g.~\cite{von1977simplicial, patriksson2015traffic}) specifically tailored to 
tackle problems with the aforementioned features. However, it is worth noting that the algorithm proposed can handle any problem of type~\eqref{prob} and can also be easily
modified in order to deal with problems having a general convex objective function. 
The reasons why we use this kind of methods are very easy to understand. Simplicial decomposition like methods, which are closely related to column generation approaches \cite{patriksson2015traffic}, are both 
well suited to deal with large-scale problems and to be used in applications where sequences of QPs need to be solved (since they can take advantage of warmstarts). 
Those tools can thus be fruitfully used in, e.g., Branch and Price like schemes for convex quadratic integer programming.

The paper is organized as follows.
In Section~\ref{sec:SD}, we describe in depth the classic simplicial decomposition framework.
In Section~\ref{sec:Strategies}, we present some strategies to improve the efficiency of the framework itself.
In Section~\ref{sec:NumericalResults}, we report our numerical experience.
Finally, in Section~\ref{sec:Conclusions}, we draw some conclusions.

\section{Simplicial Decomposition}\label{sec:SD}

Simplicial Decomposition (SD) represents a class of methods used for dealing with large scale convex problems. It was 
first introduced by Holloway in \cite{holloway1974extension} and then further studied in other papers like, e.g.,
\cite{von1977simplicial, hearn1987restricted, ventura1993restricted}. A complete overview of this kind of methods can be found in 
\cite{patriksson2015traffic}.

The method basically uses an iterative \emph{inner approximation} of the feasible set $X$.
The method can be viewed as a special case of \emph{column generation} applied to a non linear problem (we refer the reader to~\cite{desaulniers2006column} for an extensive analysis of such a method).
In practice, the feasible set $X$ is approximated with the convex hull of an ever expanding finite set $X_k= \{ \tilde{x}_1, \tilde{x}_2,\dots,\tilde{x}_m \}$ where $\tilde x_i$, $i=1,\dots,~m$ are extreme points of X.
We denote this set with $conv(X_k)$:
\begin{equation}\label{inner_descr}
conv(X_k) = \{x~|~ x=\sum_{i = 1}^m \lambda_i \tilde x_i, \sum_{i = 1}^m \lambda_i =1, \lambda_i \geq0 \}
\end{equation}

At each iteration, it is possible to add new extreme points to $X_k$ in such a way that a function reduction is guaranteed  when minimizing the objective function over the 
convex hull of the new (enlarged) set of extreme points. If the algorithm does not find at least one new point, the solution is optimal  and the algorithm terminates.

The use of the proposed method is particularly indicated when the following two conditions are satisfied:
\begin{enumerate}
 \item Minimizing a linear function over $X$ is much simpler than solving the original nonlinear problem;
 \item Minimizing the original objective function over the convex hull of a relatively small set of extreme points is 
 much simpler than solving the original nonlinear problem (i.e. tailored algorithms can be used for tackling the specific problem in our case).
\end{enumerate}
First condition is needed due to the way a new extreme point is generated. Indeed, this new point is the solution of the following linear programming problem
\begin{equation}\label{problin}
\begin{array}{ll}
\min& \nabla f(x_k)^{\top}(x-x_k)\\
\mbox{ s.t.} & x\in X\\

\end{array}
\end{equation}
where a linear approximation calculated at the last iterate $x_k$  (i.e. the solution obtained by minimizing $f$ over $conv(X_k)$~) is minimized over the original feasible set $X$.

Below, we report the detailed scheme related to the classical simplicial decomposition algorithm~\cite{bertsekas2015convex, patriksson2015traffic, von1977simplicial} (see Algorithm \ref{alg:SD}). At a generic iteration $k$ of the 
simplicial decomposition algorithm, given the set of extreme points $X_k$,
we first minimize $f$ over the set $conv(X_k)$ (Step 1), thus obtaining the new iterate $x_k$
then, at Step 2,
we generate an extreme point $\tilde x_k$ by solving the linear program \eqref{subprob}. Finally,  at Step 3, we update $X_{k}$.
\begin{algorithm}[h]
\caption{Simplicial Decomposition Algorithm}
\label{alg:SD}
  \begin{algorithmic}
  \scriptsize
  \par\vspace*{0.1cm}
  \item$\,\,\,$\hspace*{0.3truecm} \textbf{Initialization:} Choose a starting set of extreme points $X_0$.
   \par\vspace*{0.3cm}
  \item$\,\,\,$\hspace*{0.3truecm} \textbf{For} $ k=1,2,\dots $ 
   \par\vspace*{0.1cm}
  \item$\,\,\,$\hspace*{0.9truecm} \textbf{Step 1)} Generate iterate $x_k$ by solving the \textbf{master problem} 
 \begin{equation}\label{masterprob}
\begin{array}{ll}
\min& f(x)\\
\mbox{ s.t. }& x \in conv(X_k)\\
\end{array}
\end{equation}
 \item$\,\,\,$\hspace*{0.90cm} \textbf{Step 2)} Generate an extreme point $\tilde x_k$ by solving the \textbf{subproblem} 
 \begin{equation}\label{subprob}
\begin{array}{ll}
\min& \nabla f(x_k)^{\top}(x-x_k)\\
\mbox{ s.t. }& x \in X\\
\end{array}
\end{equation}
 \item$\,\,\,$\hspace*{0.9cm} \textbf{Step 3)} If $\nabla f(x_k)^{\top}(\tilde x-x_k)\geq 0$, \textbf{Stop}. Otherwise Set $X_{k+1}=X_k\cup\{\tilde x_k\}$
  \par\vspace*{0.1cm}
 \item$\,\,\,$\hspace*{0.3cm} \textbf{End For}
  \par\vspace*{0.1cm}
  \end{algorithmic}
\end{algorithm}

Finite convergence of the method is stated in the following Proposition (see, e.g., \cite{bertsekas2015convex, von1977simplicial}):
\begin{proposition}\label{finite_convergence_SD}
Simplicial Decomposition algorithm obtains a solution of Problem \eqref{prob} in a finite number of iterations.
\end{proposition}

In \cite{von1977simplicial}, a \emph{vertex dropping rule} is also used to get rid of those vertices in $X_k$ whose weight is zero in the expression of the solution $x_k$ (Step 1).
This dropping phase does not change the theoretical properties of the algorithm (finiteness still remains), but it can guarantee significant savings in terms of CPU time since it keeps the dimensions
of the master problem small.

\section{Some strategies to improve the efficiency of a simplicial decomposition framework}\label{sec:Strategies}
In this section, we discuss a few strategies that, once embedded in the simplicial decomposition framework, 
can give a significant improvement of the performances, especially when dealing with large scale quadratic problems 
with a polyhedral feasible set described by a small number of equations. 

Firstly, we present and discuss two tailored strategies to efficiently solve the master problem, 
which exploit the special structure of the generated simplices. Then, we present a couple of strategies for 
speeding up the solution of the pricing problem.

\subsection{Strategies for efficiently solving the master problem}
Here, we describe two different ways for solving the master problem. At first, we analyze an adaptive conjugate directions method
that can be used for dealing with the minimization of a quadratic function over a simplex, then we describe another tool, based 
on a \emph{ projected gradient method}, that allows us to efficiently handle 
the more general problem of minimizing a convex function over a simplex.
\subsubsection{An adaptive conjugate directions based method for solving the master}\label{ACDM}
Before describing the details related to the first method, we report a result (see e.g. \cite{pshenichnyui1978numerical}) 
for the conjugate directions method that will be useful to better understand our algorithm. 
 \begin{proposition}\label{convcg}
 Conjugate directions method makes it possible to find the minimum point of a convex quadratic function $f(x):\R^n\to \R$, and the 
 solution of the problem is obtained after less than $n$ steps.
\end{proposition}
At iteration $k$, the master problem we want to solve (Step 1 of the SD Algorithm) is the following: 
\begin{equation}\label{probmaster3}
\begin{array}{ll}
\min& f(x)~=~ x^{\top} Q x + c^{\top} x  \\
\mbox{ s.t. }&  x=\sum_{i = 1}^{k-1} \lambda_i \tilde x_i\\
&\sum_{i = 1}^{k-1} \lambda_i = 1 \\
& \lambda_i \geq 0,\\
\end{array}
\end{equation}
where the set $X_k=\{\tilde x_1,\dots, \tilde x_{k-1}\}$ represents the affine basis given by all the vertices generated in the previous iterations 
(that is we are assuming, for the sake of clarity, that all points generated so far are included in the set $X_k$: if some points,  with zero weight, have been removed with the so-called \emph{vertex dropping rule}, the method works as well).
Inspired by the approach described in \cite{von1977simplicial}, we developed a procedure that uses in an efficient way 
suitably chosen sets of conjugate directions for solving the master. The main idea is trying to reuse, as much as possible, the 
conjugate directions generated at previous iterations of the SD Algorithm.

In practice, we start from the solution of the master at iteration $k-1$, namely $x_{k-1}$,
and consider the descent direction
connecting this point with the point generated by the subproblem at iteration $k-1$, namely $\bar d_{k-1}=\tilde x_{k-1}-x_{k-1}$ (we express it 
in terms of the new coordinates of problem \eqref{probmaster3}).
Furthermore, we assume that a set of conjugate directions $D=\{d_1,\dots,d_{k-2}\}$, 
also expressed in terms of the new coordinates of problem \eqref{probmaster3},
is available from previous iterations. We then use a Gram-Schmidt like procedure to turn direction $\bar d_{k-1}$ into 
a new direction $d_{k-1}$ conjugate with respect to the set $D$.
We use the basis $B=[\tilde x_1,\dots,\tilde x_k]$ to express points $x^s=x_{k-1}$ and $x^t=x_{k-1}+d_{k-1}$
thus obtaining respectively points $\lambda^s$ and $\lambda^t$. We hence
intersect the halfline emanating from $x^s$ (and passing by $x^t$) with the boundary of the simplex in \eqref{probmaster3} by
solving the following  problem:
\begin{equation}\label{probalpha}
\begin{array}{ll}
\max& \alpha \\
\mbox{ s.t. }&  (1-\alpha)\lambda^s+\alpha\lambda^t\geq0.
\end{array}
\end{equation}
The solution of problem \eqref{probalpha} can be directly written as 
$$\alpha^*=\left(\max_i\frac{\lambda^s_i-\lambda^t_i}{\lambda^s_i}\right)^{-1}.$$
We finally define point $\lambda^p=(1-\alpha^*)\lambda^s+\alpha^*\lambda^t$ and solve the following problem 
$$\min_{\beta \in[0,1 ]}f(B[(1-\beta)\lambda^s+\beta\lambda^p]).$$
If the optimal value $\beta^* <1$ we get, by Proposition \ref{convcg}, an optimal solution for the master. Otherwise, $\beta^*=1$ and we are 
on the boundary of the simplex. In this case, we just drop those vertices  whose associated coordinates are equal to zero, 
and get a new smaller basis $B$. If $B$ is a singleton, we can stop our procedure, otherwise we minimize $f(x)$ in the new subspace
defined by $B$. In order to get  a new set of conjugate directions in the considered subspace, we use directions 
connecting point $x^*=B\lambda^*=B[(1-\beta^*)\lambda^s+\beta^*\lambda^p]$ with each vertex $\tilde x_j$ in $B$ (that is $\bar d_j=\tilde x_j-x^p$) and then use a Gram-Schmidt like procedure to make them conjugate
(we want to remark that all directions $\bar d_j$ need to be expressed in terms of the new basis $B$).
We report the algorithmic scheme below (see Algorithm \ref{alg:SDmasterACDM}).
\begin{algorithm}[h]
\caption{Adaptive Conjugate Directions based Method (ACDM)}
\label{alg:SDmasterACDM}
  \begin{algorithmic}
  \scriptsize
  \par\vspace*{0.1cm}
  \item$\,\,\,$\hspace*{0.3truecm} \textbf{Data:} Basis $B$, conjugate directions $D$, and point $x_{k-1}$
     \par\vspace*{0.1cm}
 \item$\,\,\,$\hspace*{0.30cm} \textbf{Step 1)} Set $x^s=x_{k-1}$ and $D^s=\{\bar d_{k-1}\}$
    \par\vspace*{0.1cm}
 \item$\,\,\,$\hspace*{0.30cm} \textbf{Step 2)} Select a $\bar d\in D^s$ and set $D^s=D^s\setminus\{\bar d\}$
  \par\vspace*{0.1cm}
   \item$\,\,\,$\hspace*{0.3truecm} \textbf{Step 3)} Use a Gram-Schmidt like procedure to turn 
    $\bar d$ into a conjugate direction $d^s$ with
     \item$\,\,\,$\hspace*{1.5truecm}respect to $D$
   \par\vspace*{0.1cm}
 \item$\,\,\,$\hspace*{0.30cm} \textbf{Step 4)} Express points $x^s$ and $x^t=x^s+d^s$ in terms of $B$ (that is $x^s=B\lambda^s$ and $x^t=B\lambda^t$)
  \par\vspace*{0.1cm}
 \item$\,\,\,$\hspace*{0.30cm} \textbf{Step 5)} Set 
 $$\alpha^*=\left(\max_i\frac{\lambda^s_i-\lambda^t_i}{\lambda^s_i}\right)^{-1}$$
 \item$\,\,\,$\hspace*{0.30cm} \textbf{Step 6)} Calculate point $\lambda^p=(1-\alpha^*)\lambda^s+\alpha^*\lambda^t$ and find solution $\beta^*$ of 
 the problem 
$$\min_{\beta \in[0,1 ]}f(B[(1-\beta)\lambda^s+\beta\lambda^p])$$
  \par\vspace*{0.1cm}
 \item$\,\,\,$\hspace*{0.30cm} \textbf{Step 7)} If $\beta^*<1$ then set $x^*=B[(1-\beta^*)\lambda^s+\beta^*\lambda^p]$ and $D=D\cup \{d^s\}$ go to Step 9 
 \par\vspace*{0.1cm}
 \item$\,\,\,$\hspace*{1.40cm}  Else drop vertices with $\lambda^*=0$ from $B$
 \par\vspace*{0.1cm}
 \item$\,\,\,$\hspace*{0.30cm} \textbf{Step 8)} If $B$ is a singleton then STOP 
 \par\vspace*{0.1cm}
 \item$\,\,\,$\hspace*{1.40cm}  Else set $D=\emptyset$ and for each $\tilde x_j \in B$ set $\bar d_j=\tilde x_j-x^*$ (direction represented using 
   \par\vspace*{0.1cm}
 \item$\,\,\,$\hspace*{1.40cm}  coordinates in $B$) to get a set of directions $D^s$  and go to Step 2
 \par\vspace*{0.1cm}
 \item$\,\,\,$\hspace*{0.30cm} \textbf{Step 9)} If $D^s=\emptyset$ then STOP 
 \par\vspace*{0.1cm}
 \item$\,\,\,$\hspace*{1.40cm}  Else  go to Step 2
 
  \end{algorithmic}
\end{algorithm}

Finite convergence of an SD scheme that uses Algorithm \ref{alg:SDmasterACDM} for solving the master can be obtained by using same arguments as in \cite{von1977simplicial}.
The proof is based on the fact that our polyhedral feasible set contains a finite number of simplices 
(whose vertices are extreme points of the feasible set). Since the interior of each simplex has 
at most one relative minimum and the objective function strictly decreases between two consecutive points $x_k$ and $x_{k+1}$ (keep in mind that $\nabla f(x_k)^{\top}(\tilde x^k-x^k)<0$),
no simplex can recur. Now, observing that at each iteration we get a new simplex, we have that the number of iterations must be finite.

\subsubsection{A fast gradient projection method  for solving the master}\label{FGPM}
The second approach is a Fast Gradient Projection Method (FGPM) and belongs to the family of gradient projection approaches 
(see e.g. \cite{birgin2000nonmonotone} for an overview of gradient projection approaches). 
The detailed scheme is reported below (See Algorithm \ref{alg:SDmasterFGPM}). At each iteration of the method, the new point we generate is

$$\lambda_{k+1}=\lambda_k+\beta_k(p[\lambda_k-s_k\nabla f(\lambda_k)]_\Delta-\lambda_k),$$
where $\beta_k\in(0,\rho_k]$, $\rho_k, s_k>0$ and $p[\lambda_k-s_k\nabla f(\lambda_k)]_\Delta$ is the projection 
over the master simplex in \eqref{probmaster3} of the point 
$\lambda_k-s_k\nabla f(\lambda_k)$, chosen along the antigradient. When $p[\lambda_k-s_k\nabla f(\lambda_k)]_\Delta\neq \lambda_k$, it is easy to see that the direction we get 
is a feasible descent direction.\\
The method can be used in two different ways:
\begin{itemize}
 \item[a)]we fix $s_k$ to a constant value and use a line search technique to get $\beta_k$;
 \item[b)] we fix $\beta_k$ and make a search changing $s_k$ (thus getting a curvilinear path in the feasible set).
\end{itemize}
In our algorithm we consider case $a)$ where $s_k=s>0$. 
\begin{algorithm}[h]
\caption{Fast Gradient Projection Method (FGPM)}
\label{alg:SDmasterFGPM}
  \begin{algorithmic}
  \scriptsize
  \par\vspace*{0.1cm}
  \item$\,\,\,$\hspace*{0.3truecm} \textbf{Data:} Set point $\lambda_0\in \R^{k-1}$, $\rho_0 \in [\rho_{min}, \rho_{max}]$ and a scalar value $s>0$.
   \par\vspace*{0.3cm}
  \item$\,\,\,$\hspace*{0.3truecm} \textbf{For} $k=0,1,\dots$
   \par\vspace*{0.1cm}
  \item$\,\,\,$\hspace*{0.9truecm} \textbf{Step 1)} Generate point 
 $$
\hat\lambda_k=p[\lambda_k-s\nabla f(\lambda_k)]_\Delta
$$
 \item$\,\,\,$\hspace*{0.90cm} \textbf{Step 2)} If $\hat \lambda_k =\lambda_k$  STOP; otherwise set $d_k=\hat \lambda_k-\lambda_k$
   \par\vspace*{0.3cm}
  \item$\,\,\,$\hspace*{0.90cm} \textbf{Step 3)} Choose a stepsize $\beta_k \in (0, \rho_k]$ along $d_k$  and maximum stepsize $\rho_{k+1}$ by means 
     \item$\,\,\,$\hspace*{2.0cm} of a line search 
  \par\vspace*{0.2cm}
 \item$\,\,\,$\hspace*{0.9cm} \textbf{Step 4)} Set $\lambda_{k+1}=\lambda_k+\beta_kd_k$
  \par\vspace*{0.1cm}
 \item$\,\,\,$\hspace*{0.3cm} \textbf{End For}
  \par\vspace*{0.1cm}
  \end{algorithmic}
\end{algorithm}

At each iteration, projecting the point $y_k=\lambda_k-s\nabla f(\lambda_k)$ over the simplex corresponds to solve the following problem:
$$\min_{x\in \Delta} \|x-y\|_2.$$
A fast projection over the simplex is used to 
generate the search direction \cite{condat2014fast}. This particular way of projecting a point over the simplex is basically a
 Gauss-Seidel-like variant of Michelot's variable fixing algorithm \cite{michelot1986finite}; that is, the threshold used to fix
the variables is updated after each element is read, instead of waiting for a full reading pass over the list of non-fixed elements (See \cite{condat2014fast} for further details).

A nonmonotone line search \cite{grippo1986nonmonotone} combined with a spectral steplength choice is then used at Step 3 (see \cite{birgin2000nonmonotone} for further details) to speed up 
convergence. In Algorithm \ref{NM-Arm-sc} we report the detailed scheme of the line search. Convergence of the FPGM algorithm to a minimum
follows from the theoretical results in \cite{birgin2000nonmonotone}. Therefore, the convergence of an SD method that uses FPGM to solve the master problem 
directly follows from the results in the previous sections.

    \begin{algorithm}                    
      \caption{Non-monotone Armijo line-search (with spectral steplength choice)}  
      \label{NM-Arm-sc} 
      \begin{algorithmic} 
        \par\vspace*{0.1cm}
       \scriptsize
      \item$0$\hspace*{0.5truecm} Set  $\delta \in (0,1)$, $\gamma_1\in (0,\frac 1 2)$, $M>0$ 
      \item$1$\hspace*{0.5truecm} Update $$\bar f_k = \max_{0\le i\le \min\{M,k\}} f(\lambda_{k-i})$$
      \item$2$\hspace*{0.5truecm} Set starting stepsize $\alpha =\rho_k$ and set $j=0$
      \item$3$\hspace*{0.5truecm} {\bf While} $f( \lambda_k+ \alpha d_k )> \bar f_k+\gamma_1\, \alpha\, \nabla f(\lambda_k)^{\top} d_k$
      \item$4$\hspace*{2.0truecm}set $j = j+1$ and $\alpha = \delta^j \alpha$.
      \item$5$\hspace*{0.5truecm} {\bf End While}
      \item$6$\hspace*{0.5truecm} Set $y_k=\nabla f(\lambda_k+\alpha d_k)-\nabla f(\lambda_k)$ and $b_k=\alpha d_k^{\top} y_k$
      \item$7$\hspace*{0.5truecm} If $b_k\le 0$ set $\rho_{k+1}=\rho_{max}$ else set $a_k=\alpha^2\|d_k\|^2$ and
                           $$\rho_{k+1}=min\{\rho_{max},\max\{\rho_{min},a_k/b_k\}\}$$
        \par\vspace*{0.1cm}
      \end{algorithmic}
    \end{algorithm}
  \par\smallskip\noindent

    In the FGPM Algorithm, we exploit the particular structure of the feasible set in the master, thus getting a very fast algorithm in the end. We will see 
later on that the FGPM based SD framework is even competitive with the ACDM based one, when dealing with some specific quadratic instances.


\subsection{Strategies for efficiently solving the pricing problem}
Now we describe two different strategies for speeding up the solution of the pricing problem (also called subproblem).
The first one is an \emph{early stopping} strategy that allows us to approximately solve the subproblem while guaranteeing 
finite convergence. The second one is the use of suitably generated inequalities (the so called \emph{shrinking cuts}) that both 
cut away a part of the feasible set and enable us to improve the quality of extreme points picked in the pricing phase. 
\subsubsection{Early stopping strategy for the pricing}\label{ES}
When we want to solve problem \eqref{prob} using simplicial decomposition, 
efficiently handling the subproblem is, in some cases, crucial. Indeed, 
the total number of extreme points needed to build up the final solution can be small for some real-world problem, 
hence the total time spent to solve the master problems is negligible  when compared to the total time needed to solve subproblems.
This is the reason why we may want to approximately solve subproblem \eqref{subprob} in such a way that finite convergence 
is guaranteed (a similar idea was also suggested in \cite{bertsekas2015convex}). In order to do that, we simply need to generate 
an extreme point $\tilde x_k$ satisfying the following condition:
\begin{equation}\label{es_constr} 
 \nabla f(x_k)^{\top}(\tilde x_k-x_k) \leq -\eps < 0,
\end{equation}
with $\eps>0$. Roughly speaking, we want to be sure that, at each iteration $k$, $d_k=\tilde x_k-x_k$ is a descent direction.
Below, we report the detailed scheme related to the simplicial decomposition algorithm with early stopping (see Algorithm \ref{alg:SDES}).
\begin{algorithm}[h]
\caption{Simplicial Decomposition with Early Stopping Strategy for the Subproblem}
\label{alg:SDES}
  \begin{algorithmic}
  \scriptsize
  \par\vspace*{0.1cm}
  \item$\,\,\,$\hspace*{0.3truecm} \textbf{Initialization:} Choose a starting set of extreme points $X_0$
   \par\vspace*{0.3cm}
  \item$\,\,\,$\hspace*{0.3truecm} \textbf{For} $k=0,1,\dots$
   \par\vspace*{0.1cm}
  \item$\,\,\,$\hspace*{0.9truecm} \textbf{Step 1)} Generate iterate $x_k$ by solving the \textbf{master problem} 
 $$
\begin{array}{ll}
\min& f(x)\\
\mbox{ s.t. }& x \in conv(X_k)\\
\end{array}
$$
 \item$\,\,\,$\hspace*{0.90cm} \textbf{Step 2)} Generate an extreme point $\tilde x_k\in X$ such that
 $$\nabla f(x_k)^{\top}(\tilde x_k-x_k) \leq -\eps < 0.$$
  \item$\,\,\,$\hspace*{2.00cm}   In case this is not possible, pick $\tilde x_k$ as the optimal solution of \eqref{subprob}
  \par\vspace*{0.3cm}
 \item$\,\,\,$\hspace*{0.9cm} \textbf{Step 3)}  If $\nabla f(x_k)^{\top}(\tilde x-x_k)\geq 0$, \textbf{Stop}. Otherwise set $X_{k+1}=X_k\cup\{\tilde x_k\}$
  \par\vspace*{0.1cm}
 \item$\,\,\,$\hspace*{0.3cm} \textbf{End For}
  \par\vspace*{0.1cm}
  \end{algorithmic}
\end{algorithm}
\par\smallskip\noindent

At a generic iteration $k$ we generate an extreme point $\tilde x_k$ by approximately solving the linear program \eqref{subprob}. This is done in practice by stopping the algorithm used to solve problem~\eqref{subprob} as soon as a solution satisfying constraint~\eqref{es_constr} is found. In case no solution satisfies the constraint, we simply pick the optimal solution of 
\eqref{subprob} as the new vertex to be included in the simplex at the next iteration.
\par\smallskip
Finite convergence of the method can be proved in this case as well:
\begin{proposition}\label{finite_convergence_early_stopping}
Simplicial decomposition with early stopping strategy for the subproblem  obtains a solution of Problem \eqref{prob} in a finite number of iterations.  
\end{proposition}
\par\smallskip\noindent
\textbf{Proof.} Extreme point $\tilde x_k$, obtained approximately solving subproblem \eqref{subprob}, 
can only satisfy one of the following conditions
\begin{enumerate}
 \item $\nabla f(x_k)^{\top}(\tilde x_k-x_k)\geq 0$, and subproblem \eqref{subprob} is solved to optimality. Hence we get 
$$\min_{x \in X} \nabla f(x_k)^{\top}(x-x_k)= \nabla f(x_k)^{\top}(\tilde x_k-x_k)\geq 0,$$
that is necessary and sufficient optimality conditions are satisfied and $x_k$ minimizes $f$ over the feasible set $X$;
\item $\nabla f(x_k)^{\top}(\tilde x_k-x_k)<0,$ whether the pricing problem is solved to optimality or not,
that is direction $d_k=\tilde x_k-x_k$ is descent direction and 
\begin{equation}\label{notinXk}
\tilde x_k\notin conv(X_k). 
\end{equation}
Indeed, since $x_k$ minimizes $f$ over $conv(X_k)$ it satisfies necessary and sufficient optimality conditions, that is 
$\nabla f(x_k)^{\top}(x-x_k)\geq 0$ for all $x \in conv(X_k)$.  
\end{enumerate}
From \eqref{notinXk} we thus have $\tilde x_k \notin X_k$. Since our feasible set $X$ has a finite number 
of extreme points, case 2) occurs only a finite number of times, and case 1) will eventually occur.
\hfill$\Box$

\subsubsection{Shrinking cuts}\label{Shrinking}
It is worth noticing that, at each iteration $k$, the objective function values of the subsequent 
iterates $x_{k+1}, x_{k+2},\dots,$ generated by the method will be not greater than the objective function value obtained 
in $x_k$, hence the following condition will be satisfied:
 \begin{equation}\label{condrem}
\nabla f(x_k)^{\top}(x-x_k)\leq 0. 
\end{equation}
This can be easily seen by taking into account convexity of $f$. Indeed, choosing two points $x,y\in \R^n$, we have: 
$$f(y)\geq f(x)+\nabla f(x)^{\top}(y-x).$$
Thus, if $\nabla f(x)^{\top}(y-x)>0$, we get $f(y)>f(x)$. Hence, $f(y)\leq f(x)$ implies $\nabla f(x)^{\top}(y-x)\leq0$.

We remark that all those vertices $\tilde x_i\in X_k$ not satisfying condition \eqref{condrem} have 
the related coefficient $\lambda_i=0$ in the convex combination \eqref{inner_descr} giving the master 
solution $x_k$ at iteration $k$. Proving this fact by contradiction
is easy.
Indeed, if we assume that a vertex $\tilde x_i$ is such that $\nabla f(x_k)^{\top}(\tilde x-x_k)>0 $ and the related 
$\lambda_i\neq0$, then we can build a feasible descent direction in $x_k$ thus contradicting its optimality.

We can take advantage of this property, as also briefly discussed in~\cite{bertsekas2015convex}, by adding 
the cuts described above.
The basic idea is the following: let $x_k$ be the optimal point generated by the master at a generic iteration $k$, we can hence 
add the following \emph{shrinking cut} $c_k$ to the next pricing problems:  
$$(c_k)~~~~~~~~\nabla f(x_k)^{\top}( x-x_k)\leq 0 .$$
More precisely, let 
$\{x_1, \dots, x_k\}$ be the set of optimal points generated by the master problems up to iteration $k$; then, for $k>0$, we identify as $C_k$ the polyhedron defined by all the associated shrinking cuts as follows:
$$ C_k=\{x\in {\R}^n: ~\nabla f(x_i)^{\top}( x-x_i)\leq 0, ~i=0, \dots, k-1\}.$$
(We are assuming $x_0:=\tilde{x}_0$). Therefore, at Step 2, we generate an extreme point $\tilde{x}_k$  by minimizing the linear function $\nabla f(x_k)^{\top}(x-x_k)$ over the polyhedral set $X \cap C_{k}$.
Finally, at Step 3, if $\nabla f(x_k)^{\top}(\tilde x-x_k)\geq 0$, the algorithm stops, otherwise we update $X_{k}$ by adding the point $\tilde x_k$ and $C_k$ by adding the cut $\nabla f(x_k)^{\top}( x-x_k)\leq 0$.

Below, we report the detailed scheme related to the simplicial decomposition algorithm with shrinking cuts (see Algorithm \ref{alg:SDS}).
\begin{algorithm}[h]
\caption{Simplicial Decomposition with Shrinking Cuts}
\label{alg:SDS}
  \begin{algorithmic}
  \scriptsize
  \par\vspace*{0.1cm}
  \item$\,\,\,$\hspace*{0.3truecm} \textbf{Initialization:} Choose a starting set of extreme points $X_0$
   \par\vspace*{0.3cm}
  \item$\,\,\,$\hspace*{0.3truecm} \textbf{For} $k=0,1,\dots$
   \par\vspace*{0.1cm}
  \item$\,\,\,$\hspace*{0.9truecm} \textbf{Step 1)} Generate iterate $x_k$ by solving the \textbf{master problem} 
 $$
\begin{array}{ll}
\min& f(x)\\
\mbox{ s.t. }& x \in conv(X_k)\\
\end{array}
$$
   \par\vspace*{0.3cm}
  \item$\,\,\,$\hspace*{0.90cm} \textbf{Step 2)} Generate an extreme point $\tilde x_k$ by solving the \textbf{subproblem} 
 \begin{equation}\label{subprobshrink}
\begin{array}{ll}
\min& \nabla f(x_k)^{\top}(x-x_k)\\
\mbox{ s.t. }& x \in X \cap C_k\\
\end{array}
\end{equation}
  \par\vspace*{0.2cm}
 \item$\,\,\,$\hspace*{0.9cm} \textbf{Step 3)} If $\nabla f(x_k)^{\top}(\tilde x-x_k)\geq 0$, \textbf{Stop}.  Otherwise set $X_{k+1}=X_k\cup\{\tilde x_k\}$ and
  \par\vspace*{0.1cm}
   \item$\,\,\,$\hspace*{2.0cm} set $C_{k+1}=\{x\in {\R}^n: ~\nabla f(x_i)^{\top}( x-x_i)\leq 0, ~i=0, \dots, k\}$
    \par\vspace*{0.1cm}
 \item$\,\,\,$\hspace*{0.3cm} \textbf{End For}
  \par\vspace*{0.1cm}
  \end{algorithmic}
\end{algorithm}
\par\smallskip\noindent

In practice, we implemented the algorithm with the two following variants:
\begin{itemize}
\item At the end of Step 2, after the solution of the pricing problem, we remove all shrinking cuts that are not active. In this way we are sure to have a pricing problem that is computationally tractable by keeping its size under control.
\item After a considerably large number of iterations $\bar k$, no more shrinking cuts are added to the pricing. This is done to ensure the convergence of the Algorithm.
\end{itemize}

\par\smallskip
Finite convergence of the method is stated in the following Proposition:
\begin{proposition}\label{finite_convergence_cuts}
Simplicial decomposition algorithm with shrinking cuts obtains a solution of Problem \eqref{prob} in a finite number of iterations.  
\end{proposition}
\par\smallskip\noindent
\textbf{Proof.} We first show that at each iteration the method gets a reduction of $f$ when suitable conditions are satisfied.
Since at Step 2 we get an extreme point $\tilde x_k$ by solving subproblem \eqref{subprobshrink}, if $\nabla f(x_k)^{\top}(\tilde x_k-x_k) <0$,
we have that $d_k=\tilde x_k-x_k$ is a descent direction and there exists an $\alpha_k\in (0,1]$ such that $f(x_k+\alpha_kd_k)<f(x_k)$. Since 
at iteration $k+1$, when solving the master problem, we minimize $f$ over the set $conv(X_{k+1})$ (including both  $x_k$ and $\tilde x_k$), then the minimizer $x_{k+1}$
must be such that 
$$f(x_{k+1})\leq f(x_k+\alpha_k d_k)<f(x_k).$$
Extreme point $\tilde x_k$, obtained solving subproblem \eqref{subprobshrink}, 
can only satisfy one of the following  conditions
\begin{enumerate}
 \item $\nabla f(x_k)^{\top}(\tilde x_k-x_k)\geq 0$. Hence we get 
$$\min_{x \in X \cap C_k} \nabla f(x_k)^{\top}(x-x_k)= \nabla f(x_k)^{\top}(\tilde x_k-x_k)\geq 0,$$
that is necessary and sufficient optimality conditions are satisfied and $x_k$ minimizes $f$ over the feasible set $X \cap C_k$. 
Furthermore, 
if $x \in X\setminus C_k$,
we get that there exists a cut $c_i$ with $i \in\{0,\dots,k-1\}$ such that $$\nabla f(x_i)^{\top}(x-x_i)>0.$$ Then, 
by convexity of $f$, we get 
$$f(x)\geq f(x_i)+\nabla f(x_i)^{\top}(x-x_i)>f(x_i)>f(x_k)$$
so $x_k$ minimizes f over $X$.

 \item $\nabla f(x_k)^{\top}(\tilde x_k-x_k)<0$,
that is direction $d_k=\tilde x_k-x_k$ is descent direction and 
\begin{equation}
\tilde x_k\notin conv(X_k). 
\end{equation}
Indeed, since $x_k$ minimizes $f$ over $conv(X_k)$ it satisfies necessary and sufficient optimality conditions, that is we have 
$\nabla f(x_k)^{\top}(x-x_k)\geq 0$ for all $x~\in~conv(X_k)$.
\end{enumerate}
Since from a certain iteration $\bar k$ on we do not add any further cut (notice that we can actually reduce cuts by removing the non-active ones), then 
case 2) occurs only a finite number of times. Thus case 1) will eventually occur.
\hfill$\Box$
\par\bigskip

Obviously, combining the Shrinking cuts with the Early Stopping strategy can be done (this is a part of what we actually do in practice) and finite convergence still holds for the simplicial decomposition framework.

\section{Computational results}\label{sec:NumericalResults}

\subsection{Instances description}
We used two sets of instances as test-bed: portfolio instances and generic quadratic instances. They are both described in the following subsections.

\subsubsection{Portfolio optimization problems}

We consider the formulation for portfolio optimization problems proposed by Markowitz in~\cite{markowitz}. 
The instances used have a quadratic objective function (the risk, i.e. the portfolio return variance)
and only two constraints: one giving a lower bound $\mu$ on the expected return and one representing the so call ``budget'' constraint.
The problem we want to solve is then described as follows

\begin{align}
\min_{x \in {\R}^n}\: f(x) &= x^{\top}\Sigma x\\
\mbox{s.t.} \quad r^{\top}x &\geq \mu,\nonumber\\
e^{\top}x &=1,\nonumber\\
x&\geq 0, \nonumber
\end{align}
where $\Sigma\in {\R}^{n\times n}$ is the covariance matrix, $r$ is the vector of the expected returns, 
and $e$ is the n-dimensional vector of all ones.

We used data based on time series provided in \cite{beasley}  and\cite{cesarone}. Those data are related to sets of assets 
of dimension $n = $ 226, 457, 476, 2196. The expected return and the covariance matrix are calculated by the related estimators on the time series 
related to the values of the assets.

In order to analyze the behavior of the algorithm on larger dimensional problems, we created additional instances
using data series obtained by modifying the existing ones.
More precisely, we considered the set of data with $n=2196$, and we generated bigger series
by adding additional values to the original ones: in order not to have a negligible correlation, 
we assumed that the additional data have random values close to those of the other assets. 
For each asset and for each time, we generate from 1 to 4 new values, thus obtaining 4 new 
instances whose dimensions are multiples of 2196 (that is 4392, 6588, 8784, 10980). 

For each of these 8 instances, we chose 5 different thresholds for the expected return: 
$0.006$, $0.007$, $0.008$, $0.009$, $0.01$, we thus obtained 40 portfolio optimization instances.

\subsubsection{Generic quadratic problems}

The second set of instances is of the form:
\begin{align}
\min\: f(x) &= x^{\top}Q x + c^{\top}x\\
\mbox{s. t.} \quad A x &\geq b,\nonumber\\
l \leq\, &x \leq u.\nonumber
\end{align}
with $Q\in {\R}^{n\times n}$  symmetric and positive definite matrix, $c\in {\R^n}$,  
$A\in {\R}^{m\times n}$, $b\in {\R^m}$ , $l, u\in {\R^n}$ and $-\infty<l\leq u<+\infty$. In particular, 
$Q$ was built starting from its singular value decomposition using the following procedure:
\begin{itemize}
\item the $n$ eigenvalues were chosen in such a way that they are all positive and equally distributed in the interval $(0,3]$;
\item the $n\times n$ diagonal matrix $S$, containing these eigenvalues in its diagonal, was constructed;
\item an orthogonal, $n\times n$ matrix $U$ was supplied by the QR  factorization of a randomly generated $n\times n$ square matrix;
\item finally, the desired matrix $Q$ was given by $Q=USU^{\top}$, so that it is symmetric and
its eigenvalues are exactly the ones we chose.
\end{itemize}
The coefficients of the linear part of the objective function were randomly obtained, in a small range, between $0.05$ and $0.4$, 
in order to make the solution of the problem quite sparse.

The $m$ constraints (with $m\ll n$) were generated in two different ways: step-wise sparse constraints (S) or random dense ones (R). 
In the first case, for each constraint, the coefficients associated to short overlapping sequences of consecutive variables 
were set equal to 1 and the rest equal to 0. More specifically, if $m$ is the number of constraints and 
$n$ is the number of columns, we defined $s=2*n/(m+1)$ and all the coefficients of each $i$-th constraint are zero 
except for a sequence of $s$ consecutive ones, starting at the position $1+(s/2)*(i-1)$. In the second case, 
each coefficient of the constraint matrix takes a uniformly generated random  
value in the interval $[0,1]$. The right-hand side was generated in such a way to make 
all the problems feasible: for the step-wise constraints, the right hand side was set equal to
$f*s/n$, with $0.4\leq f\leq 1$ and for a given random constraint, the corresponding 
right-hand side $b$ was a convex combination of the minimum $a_{min}$ and the maximum $a_{max}$ 
of the coefficients related to the constraint itself, that is $b=0.75*a_{min}+0.25*a_{max}$.

Each class of constraints was then possibly combined with two additional type of constraints:
a budget type constraint (b)  $e^{\top}x =1$, and a "relaxed" budget type constraints (rb) $ slb \leq e^{\top}x \leq sub$. 
Summarizing, we obtained six different classes of instances:
\begin{itemize}
\item S, instances with step-wise constraints only;
\item S-b, instances with both step-wise constraints and budget constraint;
\item S-rb, instances with both step-wise and relaxed budget constraints; 
\item R, instances with dense random constraints only;
\item R-b, instances with both dense random constraints and budget constraint;
\item R-rb, instances with both dense random and relaxed budget constraints.
\end{itemize}
For each class, we fixed $n=2000,3000,\dots,10000$, while the number of both step-wise and dense random constraints $m$ 
was chosen in two different ways:
\begin{itemize}
 \item[1)] $m=2,\,22,\,42$ for each value of $n$; 
 \item[2)] $m = n/32$, $n/16$, $n/8$, $n/4$, $n/2$ for each value of $n$.
\end{itemize}
In the first case,  we then have problems with a small number of constraints,  while, in the second case, we
have problems with a large number of constraints. 
Finally, for each class and combination of $n$ and $m$ we randomly generated five instances.
Hence, the total number of instances with a small number of constraints was 450 and the total number of  
instances with a large number of constraints was  750.

\subsection{Preliminary tests}
Here, we first describe the way we chose the  \texttt{Cplex} optimizer for solving our convex quadratic instances. 
Then, we explain how we set the parameters in the different algorithms used to solve the master problem in the SD framework. 

\subsubsection{Choice of the Cplex optimizer}
As already mentioned, we decided to benchmark our algorithm against \texttt{Cplex} version 12.6.2 (see \cite{cplex} for further details).
The optimizers that can be used in  \texttt{Cplex} for solving convex quadratic continuous problems are the following: 
primal simplex, dual simplex, network simplex, barrier, sifting and concurrent. The aim of our first test was to identify, 
among the 6 different options, which is the most efficient for solving instances with a dense $Q$ and $n \gg m$. 

In Table~\ref{tab:cplexcomparison}, we present the results concerning instances with 42 constraints and 
three different dimensions $n$: 2000, 4000 and 6000. We chose problems with a small number of constraints  
in order to be sure to pick the best \texttt{Cplex} optimizer for those problems where the SD framework is 
supposed to give very good performances. For a fixed $n$, three different instances 
were solved of all six problem types. So, each entry of Table~\ref{tab:cplexcomparison} represents the averages 
computing times over 18 instances. 
A time limit of 1000 seconds was imposed and in brackets we report (if any) the number of instances 
that reached the time limit.

\input{cplex_comparison.tex}

The table clearly shows that the default optimizer, the barrier and the concurrent methods give poor performances 
when dealing with the quadratic programs we previously described. 
On the other side, the simplex type algorithms and the sifting algorithm seem to be very fast for those instances.
In particular, sifting gives the overall best performance.
Taking into account these results, we decided to use the \texttt{Cplex} sifting optimizer as the baseline method in our experiments.
It is worth noticing that the sifting algorithm is specifically conceived by \texttt{Cplex} to   
deal with problems with $n \gg m$, representing an additional reason for comparing our algorithmic framework against
this specific \texttt{Cplex} optimizer.

\subsubsection{Tolerance setting when solving the master problem}\label{tolerance}

We have three options available for solving the master problem in the SD framework: ACDM, FGPM and \texttt{Cplex}. 
In order to identify the best choice, we need to properly set tolerances for those methods.
When using \texttt{Cplex} as the master solver, we decided to keep the tolerance to its default value (that is $1E10-6$).
The peculiar aspect of ACDM is that no tolerance needs to be fixed a priori.
On the other hand, with FGPM, the tolerance setting phase is very importance since, as we will see, 
it can significantly change the performance of the algorithm in the end.

In Table~\ref{tab:errFGPM}, we 
compare the different behaviors of our SD framework for the
three different choices of master solver. Each line of the table represents the average values concerning the 54 instances 
used in the previous experiment.
Column ``T'' represents the time (in seconds) spent by the algorithms. ``Er'' and ``Max Er''  represent the average and maximum relative errors 
with respect to the value found by \texttt{Cplex} (using sifting optimizer). ``Ei'' and ``Max Ei'' represent the average and maximum distance (calculated using $\ell_{\infty}$ norm) 
from the solution found by \texttt{Cplex}. 
In the last column, ``Dim'' represents the dimension of the final master program.

\input{errFGPM.tex}

By taking a look at the table, we can easily see that the ACDM based SD framework gets the best results in terms of errors with respect to
\texttt{Cplex}. We can also see that the performance of the FGPM based one really changes depending on the tolerance chosen. 
If we want to get for FGPM the same errors as ACDM, we need to set the tolerance to very low values, thus considerably slowing down the 
algorithm. In the end, we decided to use a tolerance of $10E-6$ for FGPM, 
which gives a good trade-off between computational time and accuracy. This means anyway that we gave up precision to keep the algorithm fast 
with respect to ACDM.

\subsection{Numerical results related to the complete testbed}
Now, we analyze the performances of our SD framework when choosing different options for 
both solving the master and the pricing problem. 
Summarizing, we used three different methods for solving the master: \texttt{Cplex}, ACDM and FGPM. 
As for the pricing problems, we considered the use of the Early stopping (E) technique described in Section~\ref{ES} 
and the Shrinking cuts (Cuts), described in Section~\ref{Shrinking}. 
A further option we used in the SD framework is the use of the LP sifting optimizer in the solution of the pricing problem
instead of the standard one. We also notice that, in order to save CPU time, 
the vertex dropping rule described in Section \ref{sec:SD} was always used in the framework.
Summing up, for each choice of the master solver, we compared 8 different sets of options related to the pricing solver, indicated as follows:\\
\begin{tabular}{ll}
$\bullet$ Default (D) & $\bullet$ Sifting (Sif)\\
$\bullet$ Cuts (C) & $\bullet$ Sifting+Cuts (Sif-C)\\
$\bullet$ Early stopping (E) & $\bullet$ Sifting+Early stopping (Sif-E)\\
$\bullet$ Cuts+Early stopping (CE) & $\bullet$ Sifting+Cuts+Early stopping (Sif-CE)
\end{tabular}

\subsubsection{Portfolio optimization instances}
\begin{figure}[htp]
\includegraphics[height=0.24\textheight]{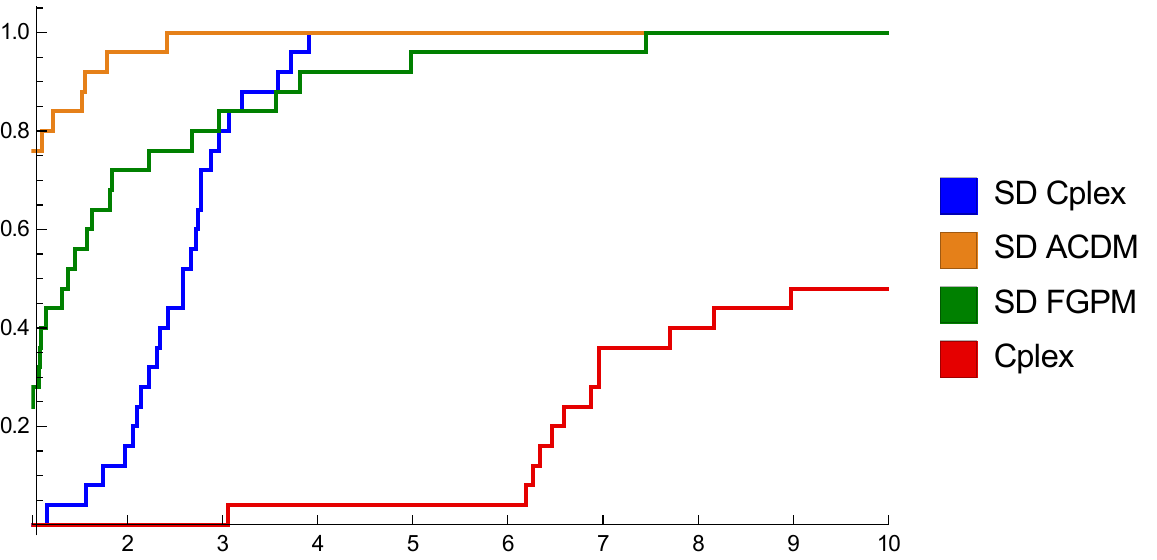}
\caption{Performance profiles for Portfolio instances - master solvers.}
\label{port_all}
\end{figure}

Firstly, we tested the framework on the portfolio optimization instances. For the largest values of $n$, that is $n>2000$, 
in Figure~\ref{port_all}, we show the performance profiles related to the framework with the three different master solvers 
and the sifting \texttt{Cplex} optimizer. We produced the performance profiles according to \cite{Dolan2002} and using the software \emph{Mathematica} version 10.2 (see \cite{mathematica} for further details). For each SD method, the best pricing option was considered. We can easily notice that
the best choice in terms of master solver is ACDM.

\begin{figure}[htp]
\includegraphics[height=0.24\textheight]{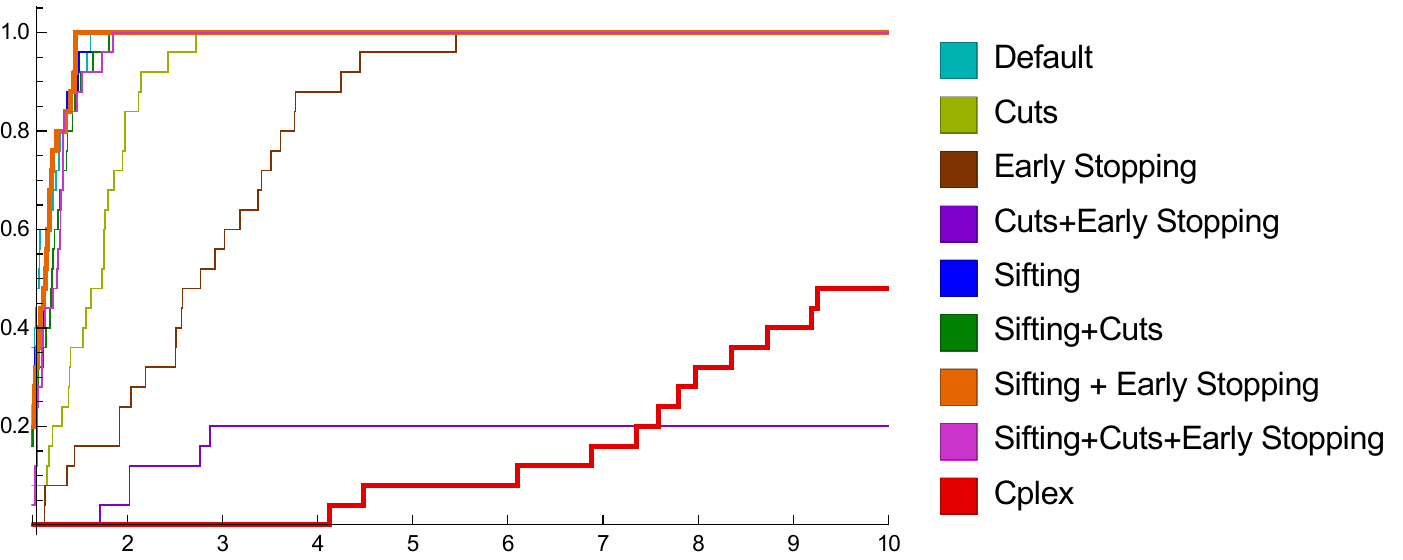}
\caption{Performance profiles for Portfolio instances - pricing options (SD ACDM).}
\label{port_ACDM}
\end{figure}

In Figure~\ref{port_ACDM}, we report the performance profiles related to the best master solver (that is ACDM)
for the different pricing options. As we can see, the best option for the pricing is Sifting + Early Stopping.

\subsubsection{Generic quadratic instances}
\begin{paragraph}{Small number of constraints (GS)}
We first analyze the results related to the generic quadratic problems with a small number of constraints. 
The performance profiles reported in Figures  ~\ref{lowm_all} and~\ref{lowm_ACDM} are related to all classes of constraints described before.

\begin{figure}[htp]
\includegraphics[height=0.24\textheight]{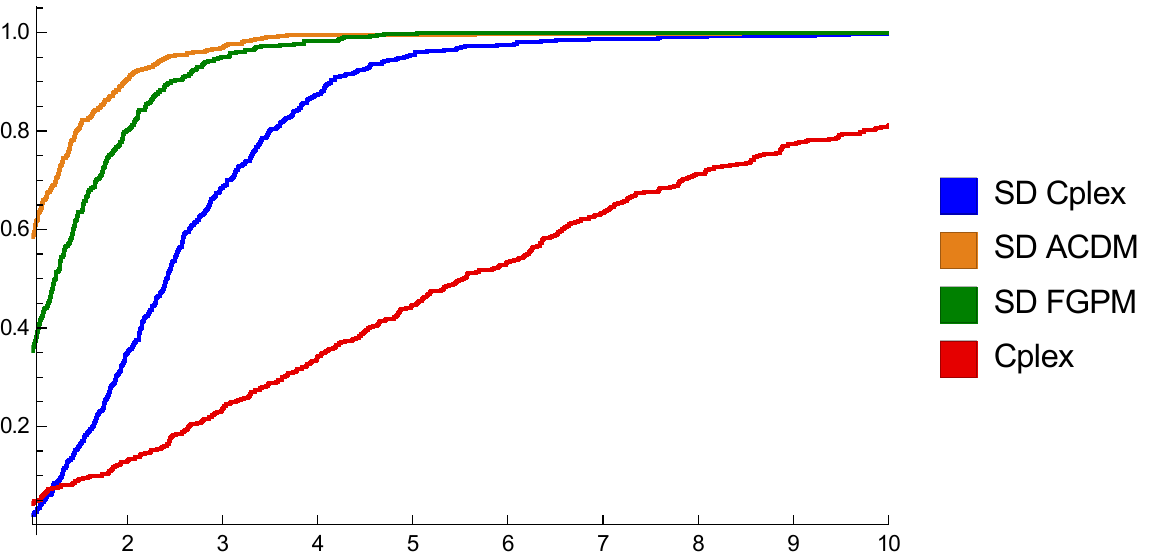}
\caption{Performance profile for GS instances - master solvers.}
\label{lowm_all}
\end{figure}

\begin{figure}[htp]
\includegraphics[height=0.24\textheight]{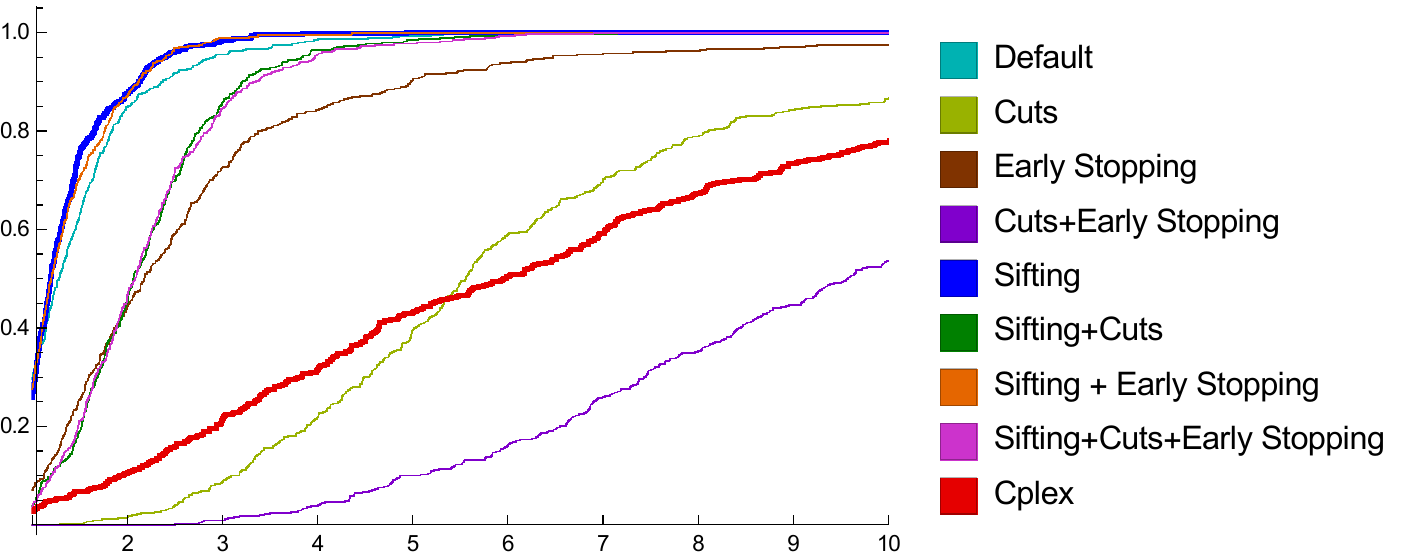}
\caption{Performance profile for GS instances - pricing options (SD ACDM).}
\label{lowm_ACDM}
\end{figure}

Similarly to the results on portfolio instances, ACDM represents the best choice for solving the master; 
sifting is the best pricing option for ACDM, whose performance does not change a lot with the addition of the Early stopping strategy.
\end{paragraph}

\begin{paragraph}{Large number of constraints (GL)} Here, we analyze the results obtained for instances with a 
larger number of constraints. We need to keep in mind, anyway, that our SD framework works well only when the number 
of constraints is significantly smaller than the number of variables.

Analogously as before, in Figures~\ref{largem_all} and ~\ref{largem_FPGM}, 
we analyze the performances of the framework for the different choices of 
master solvers and pricing options. The results are compared using performance profiles (we consider all types of constraints in the analysis).

\begin{figure}[htp]
\includegraphics[height=0.24\textheight]{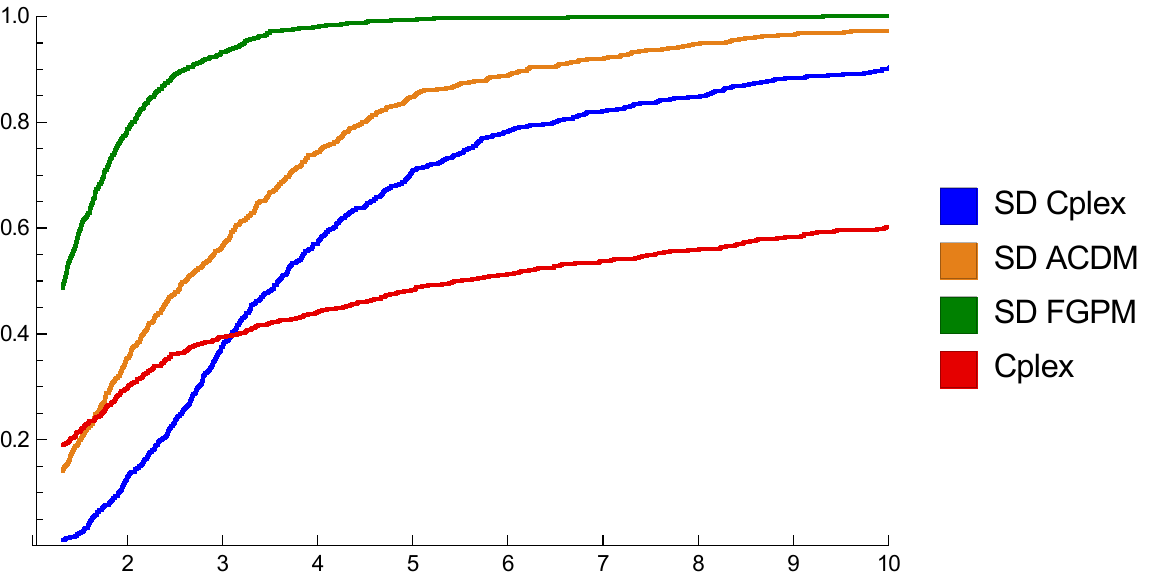}
\caption{Performance profile for GL instances - master solvers.}
\label{largem_all}
\end{figure}

\begin{figure}[htp]
\includegraphics[height=0.24\textheight]{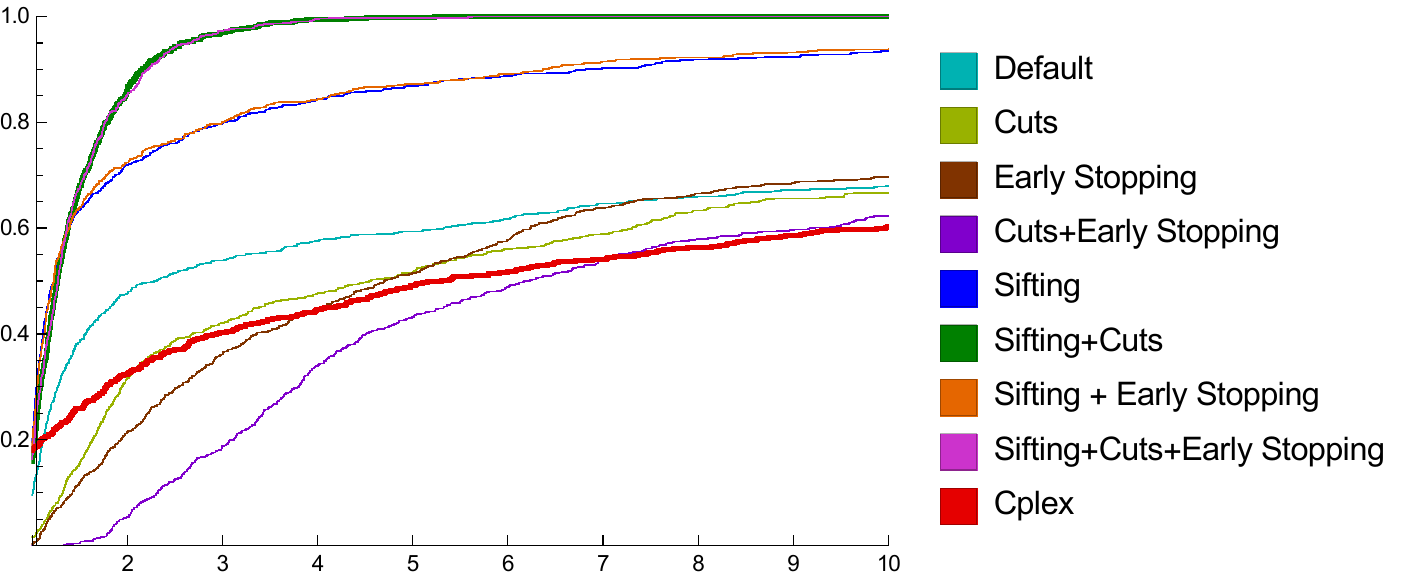}
\caption{Performance profile for GL instances - pricing options (SD FGPM).}
\label{largem_FPGM}
\end{figure}

In this case something different happens:
\begin{itemize}
\item \texttt{Cplex}, SD ACDM  and SD Cplex, reach the time limit on some instances;
\item the best master solver now is clearly FGPM. SD  ACDM and SD  \texttt{Cplex} are even worse than \texttt{Cplex} 
in terms of efficiency, but are more robust.
\item the best pricing option is Sifting + Cuts. So, the use of this type of cuts, which was ineffective when dealing with portfolio and GS instances, 
significantly improves the performance here. We get the same improvement when using the other master solvers too.
\end{itemize}


As we said, some of the algorithms did not solve all the instances within the time limit of 1000 seconds.
In particular, 122 instances out of 750 were not solved by \texttt{Cplex},  only 8 by SD \texttt{Cplex} and 
4 by SD ACDM (both with the Sifting + Cuts option); SD FGPM (with the Sifting + Cuts option), on the other side, solved all the instances within the time limit.

The biggest difference with respect to the previous results is that now the performance profiles 
of SD FGPM are better than those obtained using the other methods (and, in particular, better than SD ACDM, which was the best one so far). 
The reason why this happens will become clearer later on (see Section~\ref{deeperanalysis}), 
when an in-depth analysis of the results will be shown. We anyway need to keep in mind that SD ACDM usually gives 
better solutions in terms of errors when compared with SD FGPM.

Even though the overall results clearly say that the SD method is better than the \texttt{Cplex} sifting optimizer, 
for a subset of instances this is not always true. 
Indeed, after analyzing  the results for each specific type of instance, we noticed that for the GL problems with dense random constraints 
and with the addition of the budget constraint, the results of our SD framework were a bit worse than those obtained using \texttt{Cplex}. 
Performance profiles related to those results are reported Figure~\ref{qp5_largem_all}.

\begin{figure}[htp]
\includegraphics[height=0.24\textheight]{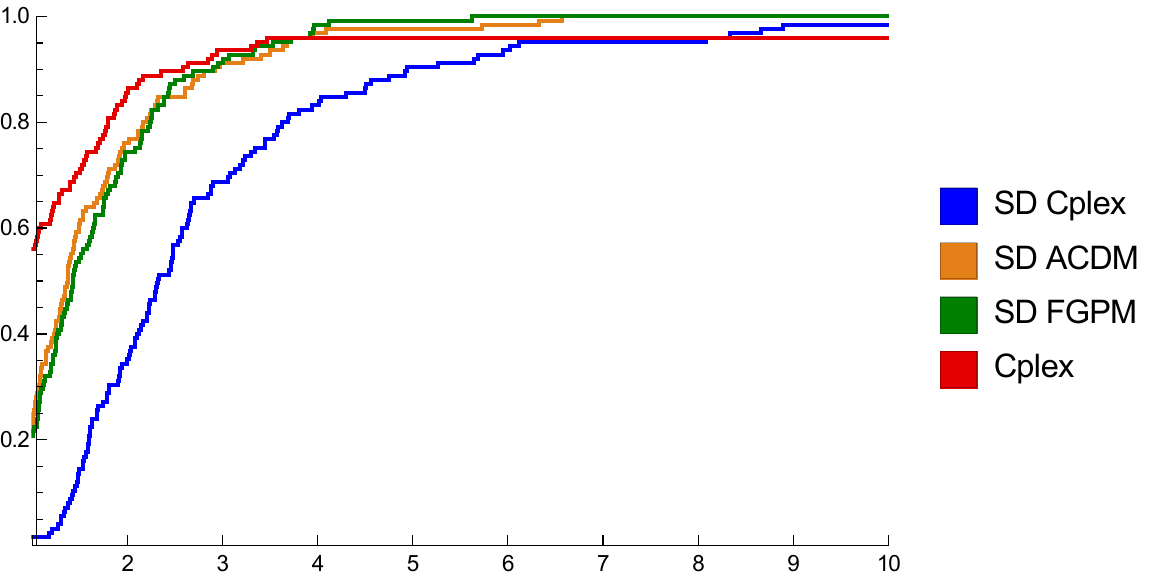}
\caption{Performance profile for GL instances (Rb constraints).}
\label{qp5_largem_all}
\end{figure}

\end{paragraph}

\subsection{CPU time usage in the SD framework}
Now we analyze the way CPU time is used in the SD framework, that is we show the average CPU time needed for preprocessing data,
solving the master problems and solving the pricing problems (failures are not considered in the analysis). In Figure~\ref{fig:timerep}, we report the aggregated results
over all the solved instances. In each figure, we report the time spent by SD in the preprocessing phase of the algorithm (\emph{preprocessing}), 
in the solution of the master and pricing problem. The solving time of both the pricing and master problem is split in the time needed 
to update the data structure (\emph{updating}) and the time needed to solve to problem (\emph{solvers}). 
For each figure we provide also the average computing time over the whole testbed.
\input{piecharts.tex}
Table~\ref{tab:num_its} contains the average number of iterations needed to achieve 
the solution and the average dimension of the last master program, which is the dimension of the optimal face in the original domain.
\input{num_its.tex}
In particular, one can note that the dimension of the last master is greater when using SD Cplex (i.e., both  SD ACDM and 
SD FGPM tend to find sparser solution than SD Cplex). From these tables it is also clear that SD FGPM saves time with respect to SD ACDM 
mainly because it needs less iterations to satisfy the stopping criterion.
This is mainly due to the tolerance chosen for FGPM. As we have seen before, the choice made
gives, on one side, better results in terms of CPU time, but, on the other side, it gives a slight deterioration  
of the solution quality in all instances. It is also worth noticing that:
\begin{itemize}
 \item the reduction obtained by SD FGPM with respect to
SD ACDM, in terms of average number of iterations,  is around $25\%$ (see Table~\ref{tab:num_its});
\item the average CPU time needed for solving the pricing in SD FGPM is nearly halved with respect to SD ACDM.
\end{itemize}
The  CPU time reduction in the pricing phase seems to be caused by two different factors: pricing problems in SD FGPM get very similar very soon
thus requiring less iterations of the LP method to be solved;  shrinking cuts are more effective when embedded into SD FGPM.

Summarizing, SD ACDM is more precise and reliable. Furthermore, it gives the best performances when dealing with portfolio an GS instances.
SD FGPM seems to be a good choice when dealing with GL instances (we find good solutions in a smaller amount of time).

\subsection{In-depth analysis}\label{deeperanalysis}
\begin{figure}[htp]
\includegraphics[width=0.6\textwidth]{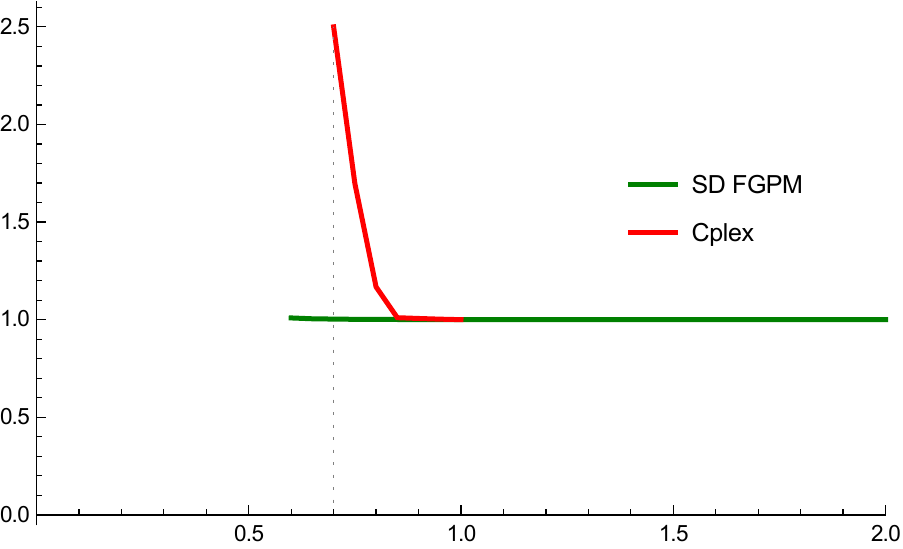}
\caption{Objective function decay - Objective function ratio (y-axis) and CPU time ratio (x-axis) - SD FGPM vs Cplex.}
\label{av_glob}
\end{figure}

\begin{figure}[htp]
\includegraphics[width=0.6\textwidth]{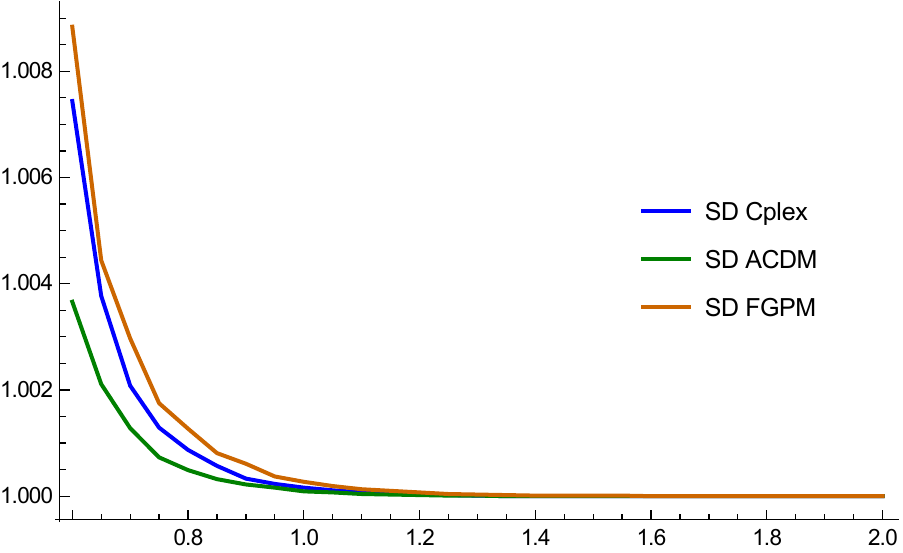}
\caption{Objective function decay - Objective function ratio (y-axis) and CPU time ratio (x-axis) - SD solvers comparison.}
\label{avSDs}
\end{figure}

\begin{figure}[htp]
\includegraphics[width=0.6\textwidth]{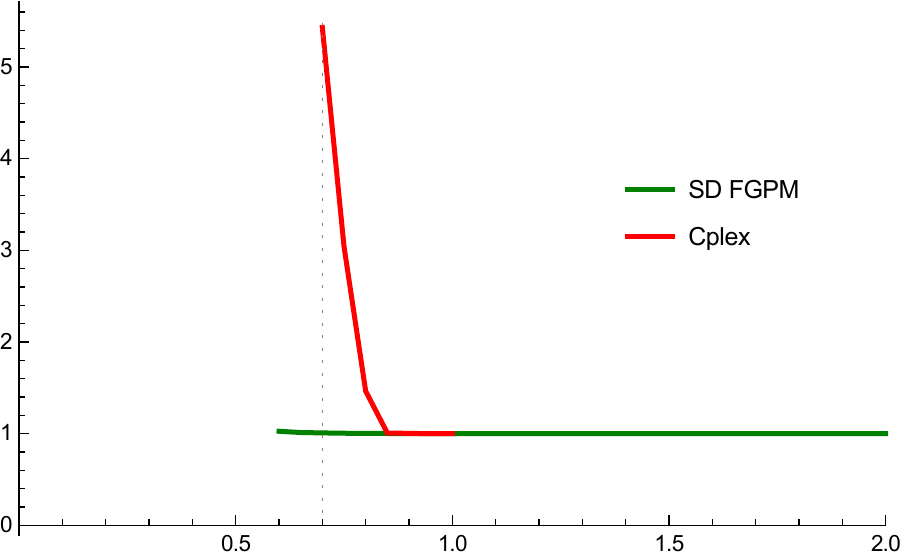}
\caption{Objective function decay - Objective function ratio (y-axis) and CPU time ratio (x-axis) - GL instances (Rb constraints).}
\label{av_5}
\end{figure}

In order to better analyze the behavior of the SD framework, we show now how the objective function value changes  
with respect to the elapsed time. 
Since we want to get meaningful results, we only consider those instances solved 
in more than 10 seconds (but always within the time limit of 1000 seconds).
In particular, we consider instances with random dense constraints and we take a set of 25 instances for each of the three types 
of additional constraints.
Hence, we plot 
\begin{itemize}
 \item on the \emph{x-axis} the \textbf{CPU time ratio}, that is the CPU time elapsed  divided by the overall time needed by   
\texttt{Cplex} to get a solution on the same instance.
\item on the \emph{y-axis} the \textbf{objective function ratio}, that is the objective function value  divided by the optimal value 
obtained by \texttt{Cplex} on the same instance.
\end{itemize}
All the results are averaged over the whole set of instances. For the SD framework, we plot the results up to twice the time needed by \texttt{Cplex}
to get a solution.
In the analysis, we always consider the setting "Sif-CE", which considers all the pricing options (and gives same performance as the best one).
Figures~\ref{av_glob} and~\ref{avSDs} show the overall results for the 75 instances considered: the first figure shows the comparison between 
\texttt{Cplex} and SD FGPM, while the second one shows the comparison of the three different SD framework versions.
From the comparison of \texttt{Cplex} and SD FGPM, it is easy to notice that SD gets a good objective function value very soon. Indeed, at a CPU time
ratio 0.6 (i.e., $60\%$ of the overall \texttt{Cplex} CPU time) corresponds an objective function ratio slightly bigger than 1 for SD FGPM, while 
at the same CPU time ratio \texttt{Cplex} still needs to find a feasible solution. \texttt{Cplex} gets a first feasible solution for a CPU time
ratio equal to 0.7  (in this case the objective function ratio is bigger than 2.5), and it obtains an objective function ratio close to 1 only 
for a CPU time ratio bigger than 0.8. By taking a look at the comparison of the three different versions of our SD framework, 
we notice that SD FGPM actually takes longer than the others to get an objective function ratio close to 1. The better results 
obtained for SD FGPM hence depend, as we already noticed, on the way we choose the tolerance in the master solvers. 
Finally, in Figure~\ref{av_5}, we report the plots related to those instances where \texttt{Cplex} outperforms the SD framework.
Once again, we can see that SD FGPM gets a good objective function ratio very soon, while \texttt{Cplex} takes much longer to obtain
a similar ratio.

\vspace{2mm}

\section{Conclusions}\label{sec:Conclusions}

We presented an efficient SD framework to solve continuous convex quadratic problems. 
It embeds two ad-hoc methods for solving the master problem, namely an
adaptive conjugate directions based method and a fast gradient projection method.
Furthermore, three different strategies to speed up the pricing are included:  
an early stopping technique, a  method to shrink the feasible region based on some specific 
cuts, and a sifting strategy to solve the pricing problem.

We showed, through a wide  numerical experience,  that our algorithm is better than \texttt{Cplex}
when dealing with instances with a dense Hessian matrix and with a number of constraints considerably  smaller than the number of variables.

In Table~\ref{overview}, we summarize the recommended settings with respect to the instances we solved. 
For portfolio instances (\textbf{PORTFOLIO}), the best master solver is ACDM and the best pricing option is sifting with early stopping.
For generic quadratic instances with a small number of constraints (\textbf{SMALL m}) the best master optimizer is again ACDM
and the the best pricing option is sifting. Finally, for generic quadratic instances with a large number of constraints (\textbf{LARGE m})
the cuts play an important role (best pricing option is sifting with cuts) and FGPM is the best master solver. 

\input{overview.tex}

\bibliographystyle{spmpsci} 
\bibliography{paper_SD_QP}

\end{document}

%% file: cplex_comparison.tex
\begin{table}[htb]
\centering

\begin{tabular}{S[table-format=4.1]S[table-format=4.2]S[table-format=4.2]S[table-format=2.1]S[table-format=2.1]S[table-format=4.2]S[table-format=4.2]S[table-format=4.2]}

\toprule

$n$	& $\text{Default}$		&	$\text{Primal}$	&	$\text{Dual}$	&	$\text{Network}$	&	$\text{Barrier}$	&	$\text{Sifting}$	&	$\text{Concurrent}$	\\
\midrule
2000	&	72.2			&	1.6	&	1.6	&	1.6	&	84.2			&	2.0	&	89.0		\\
4000	&	641.8\hspace{0.01cm}(2)	&	12.7	&	13.9	&	13.9	&	618.0\hspace{0.01cm}(2)	&	11.5	&	689.4\hspace{0.01cm}(2)	\\
6000	&	1000.0\hspace{0.01cm}(18)	&	31.5	&	30.7	&	30.5	&	1000.0\hspace{0.01cm}(18)	&	26.3	&	1000.0\hspace{0.01cm}(18)	\\

\bottomrule
\end{tabular}
\caption{Comparison among the different Cplex optimizers}
\label{tab:cplexcomparison}
\end{table}

%% file: errFGPM.tex
\begin{table}[htb]
\centering
\begin{tabular}{c*{7}{r}}
\toprule

Solver	&Tol	&	T (s)	&	Er	&	Max Er	&	Ei	&	Max Ei	&	Dim	\\
\midrule
\multirow{5}*{SD FGPM}	&	1E-02	&	0.25	&	8.64E-02	&	2.67E-01	&	2.24E-02	&		5.04E-02	&	9.9	\\
			&	1E-04	&	1.15	&	2.21E-04	&	6.79E-04	&	7.80E-04	&	1.44E-03	&	55.6	\\
			&	\textbf{1E-06}	&	2.46	&	5.65E-07	&	2.63E-06	&	5.72E-05	&	1.86E-04	&	102.2	\\
			&	1E-08	&	6.09	&	5.98E-09	&	1.15E-07	&	4.61E-06	&	1.88E-05	&	114.0	\\
			&	1E-10	&	9.81	&	2.35E-09	&	4.59E-08	&	3.48E-06	&	2.16E-05	&	113.4	\\
SD Cplex		&	1E-06	&	4.66	&	8.86E-09	&	4.26E-08	&	5.50E-06	&	2.46E-05	&	156.0	\\
SD ACDM		&	None		&	3.63	&	1.53E-09	&	1.97E-08	&	2.65E-06	&	1.99E-05	&	113.1	\\
Cplex		&			&	4.29&			&			&			&			&		\\
\bottomrule
\end{tabular}
\caption{Comparison for the three different choices of master solver (with Cplex we indicate the results obtained with sifting optimizer).}
\label{tab:errFGPM}

\end{table}

%% file: piecharts.tex
\begin{figure}[htp]
\centering
\subfigure[SD Cplex. Average CPU Time = 27.5 s]
{\includegraphics[width=0.7\textwidth]{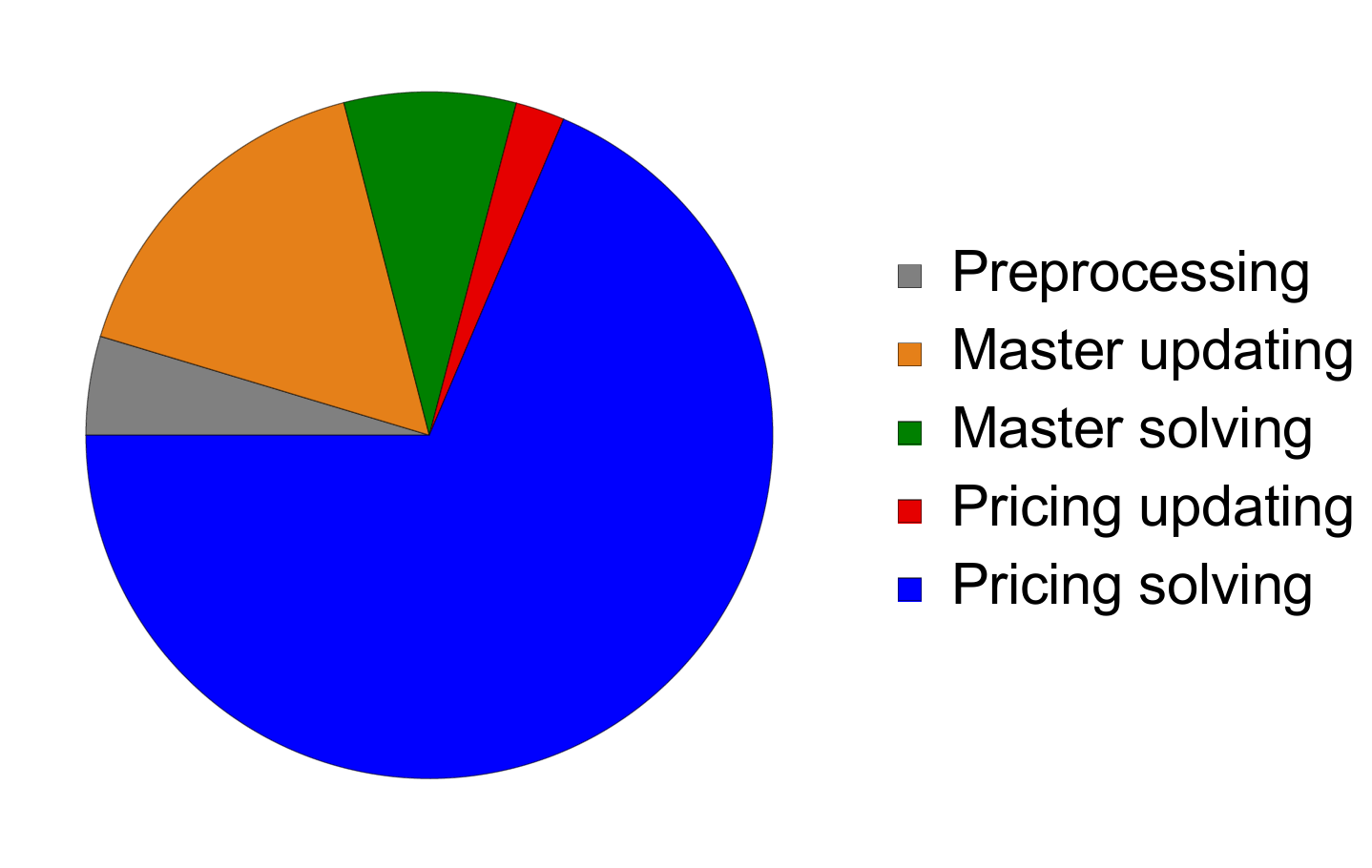}}
\hspace{5mm}
\subfigure[SD ACDM. Average CPU Time = 19.6 s]
{\includegraphics[width=0.7\textwidth]{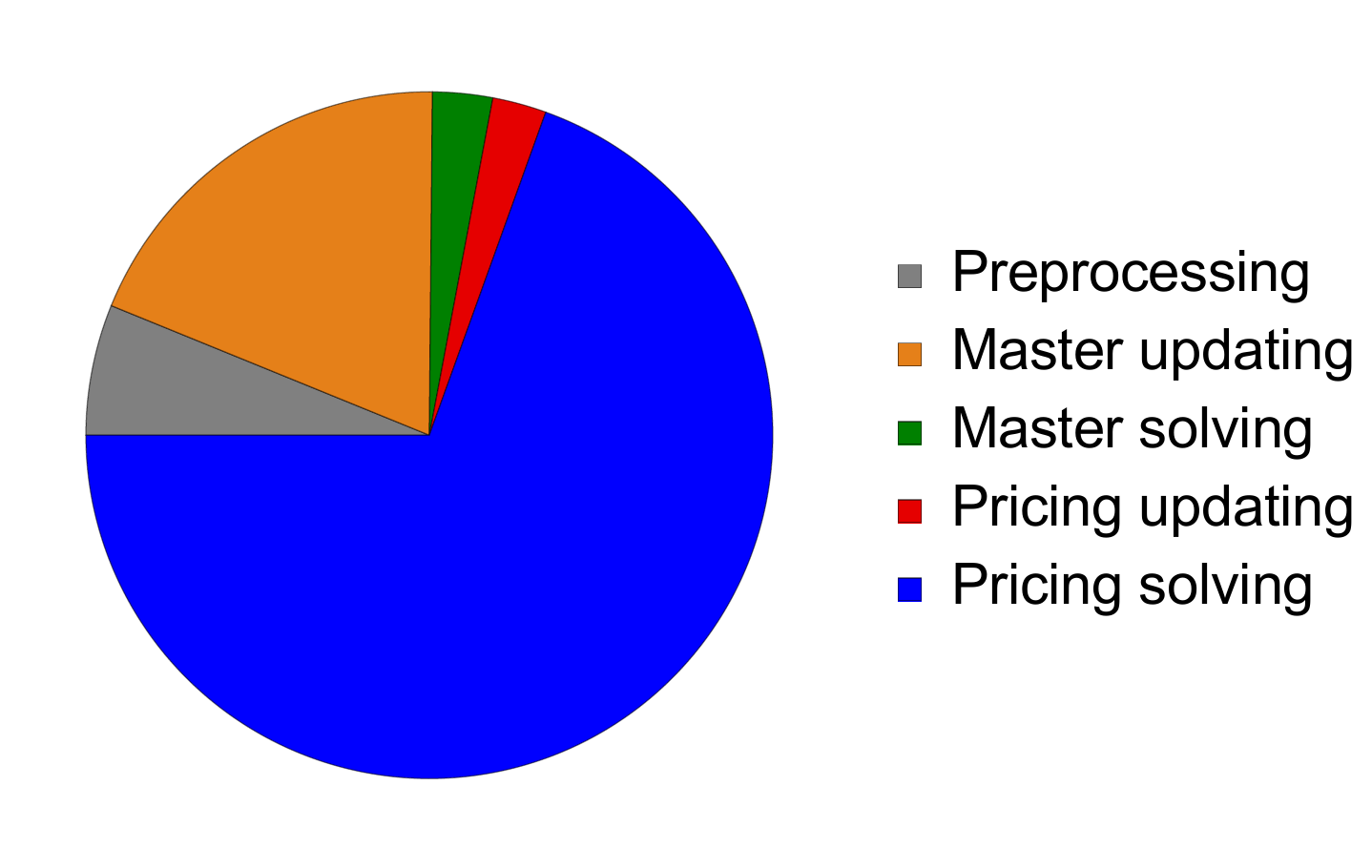}}
\hspace{5mm}
\subfigure[SD FGPM. Average CPU Time = 11.3 s]
{\includegraphics[width=0.7\textwidth]{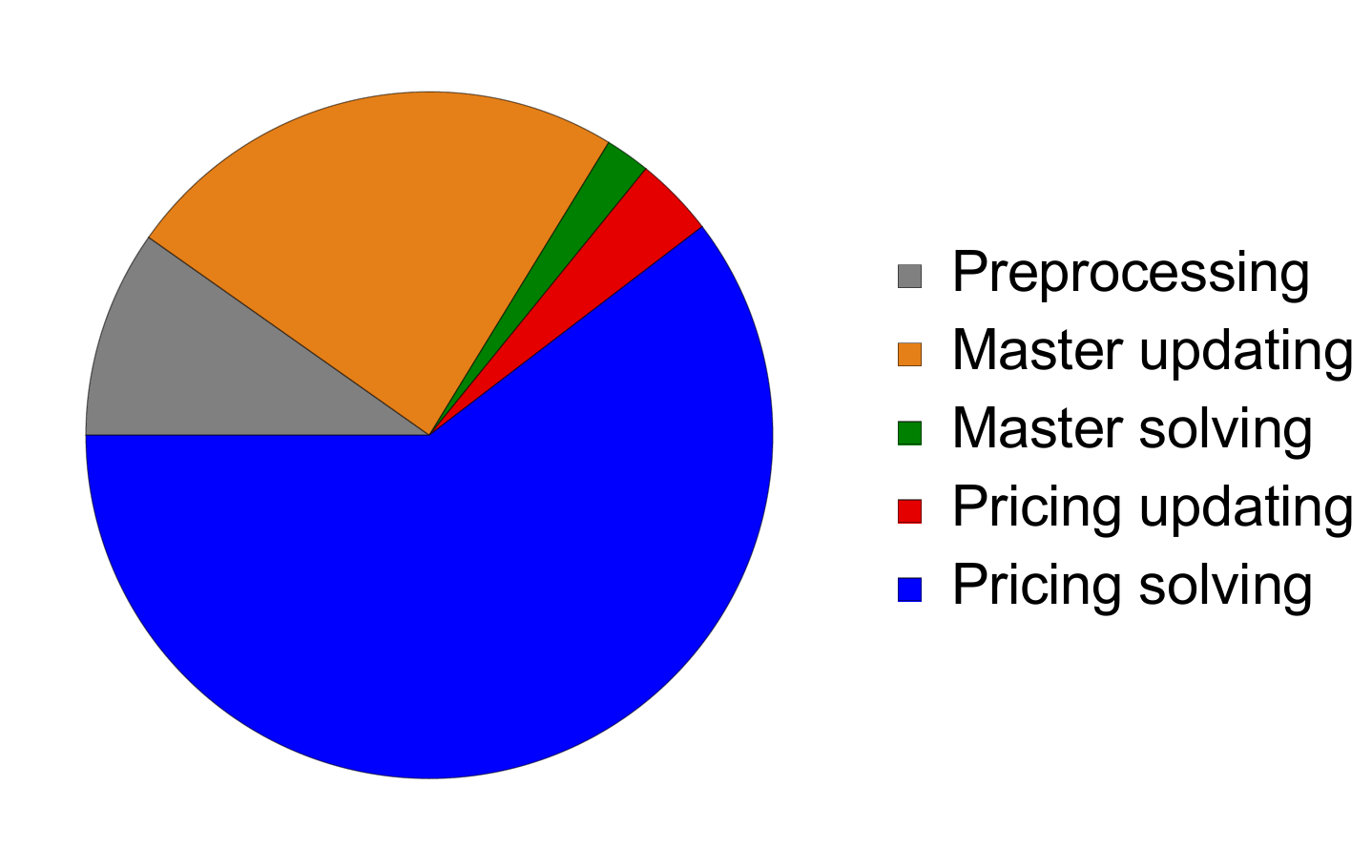}}
\hspace{5mm}
\caption{CPU time pie charts (Cplex Average CPU time = 39.8 s).}
\label{fig:timerep}
\end{figure}

%% file: num_its.tex
\begin{table}[htb]
\centering
\begin{tabular}{*{4}{r}}
\toprule
	&	\multicolumn{3}{c}{Portfolio Instances}\\
\midrule
	&	SD Cplex	&	SD ACDM	&	SD FGPM	\\
\midrule
iterations 	&	108.3&	111.2	&	99.4		\\
final dimension	&	109.3	&	79.8	&	77.7		\\
\midrule
	&	\multicolumn{3}{c}{GS Instances}	\\
\midrule
	&	SD Cplex	&	SD ACDM	&	SD FGPM	\\
\midrule
iterations 	&	171.4	&	172.4	&	145.7	\\
final dimension	&	168.3	&	137.4	&	130.7	\\
\midrule
	&	\multicolumn{3}{c}{GL Instances}	\\
\midrule
	&	SD Cplex	&	SD ACDM	&	SD FGPM	\\
\midrule
iterations 	&	132.9	&	118.4	&	89.3	\\
final dimension	&	126.6	&	85.5		&	70.7		\\
\bottomrule
\end{tabular}

\caption{Comparison for the three choices of master solvers. Average number of iterations and average dimension of the last master.}
\label{tab:num_its}

\end{table}

%% file: overview.tex
\begin{table}[htb]
\centering
\begin{adjustbox}{max width=\textwidth}
\begin{tabular}{ c r c c c c c c }
\toprule
\multicolumn{2}{c}{}&\multicolumn{2}{c}{Master solvers}&\phantom{}&\multicolumn{3}{c}{Pricing options}\\
\cmidrule{3-4} \cmidrule{6-8}
\multicolumn{2}{c}{}&ACDM&FGPM&&SIFTING &EARLY ST. & CUTS \\
\midrule
&\textbf{PORTFOLIO} & {\LARGE\checkmark} & {\LARGE$\times$}& &{\LARGE\checkmark} & {\LARGE\checkmark} & {\LARGE$\times$} \\
Instances& \textbf{SMALL m} & {\LARGE\checkmark} & {\LARGE$\times$}&& {\LARGE\checkmark} & {\LARGE$\times$} &{\LARGE$\times$} \\
& \textbf{LARGE m} & {\LARGE$\times$} & {\LARGE\checkmark}&& {\LARGE\checkmark} & {\LARGE$\times$} & {\LARGE\checkmark} \\
\bottomrule
\end{tabular}
\end{adjustbox}
\caption{Best settings overview.}
 \label{overview}
\end{table}